\documentclass[12pt]{article}
\usepackage{latexsym,amsfonts,amssymb}
\usepackage[dvips]{color}
\usepackage{colordvi,multicol}

\def \cal{\mathcal}

\newtheorem{thm}{Theorem}[section]
\newtheorem{cor}[thm]{Corollary}
\newtheorem{lem}[thm]{Lemma}
\newtheorem{pro}[thm]{Proposition}
\newtheorem{defi}[thm]{Definition}
\newtheorem{rem}[thm]{Remark}
\newtheorem{exa}[thm]{Example}

\date{}
\begin{document}
\title{\bf Hunt's Hypothesis  (H) for Markov Processes: Survey and Beyond}
\author{}
%Ze-Chun Hu, Wei Sun and Li-Fei Wang

\maketitle

 \centerline{Ze-Chun Hu} \centerline{\small College
of Mathematics}
\centerline{\small Sichuan University}
\centerline{Chengdu, 610065, China}
\centerline{\small E-mail: zchu@scu.edu.cn}

\vskip 1cm \centerline{Wei Sun} \centerline{\small Department of
Mathematics and Statistics}
\centerline{\small Concordia University}
\centerline{\small Montreal, H3G 1M8, Canada} \centerline{\small
E-mail: wei.sun@concordia.ca}

\vskip 1.5cm
\begin{abstract}

\noindent The goal of this paper is threefold. First, we survey
the existing results on Hunt's hypothesis (H) for  Markov
processes and Getoor's conjecture for L\'{e}vy processes. Second,
we investigate (H) for multidimensional L\'{e}vy processes from
the viewpoints of projections and energy, respectively. Third, we
present a few open questions for further study.

\end{abstract}

%\tableofcontents

\section{Introduction}

Let $E$ be a locally compact space with a countable base and $X=(X_t, P^x)$ be a standard Markov process on $E$ as
described in Blumenthal and Getoor \cite{BG68}. Denote by ${\cal
B}$ and ${\cal B}^n$ the family of all Borel measurable subsets
and nearly Borel measurable subsets of $E$, respectively. For
$D\subset E$, we define the first hitting time of $D$ by
$$
T_D:=\inf\{t>0:X_t\in D\}.
$$

A set $D\subset E$ is called {\it thin} if there exists a set $C\in
{\cal B}^n$ such that $D\subset C$ and {$P^x(T_C=0)=0$} for any
$x\in E$. $D$ is called {\it semipolar} if
$D\subset\bigcup_{n=1}^{\infty}D_n$ for some thin sets
$\{D_n\}_{n=1}^{\infty}$. $D$ is called {\it polar} if there exists a
set $C\in {\cal B}^n$ such that $D\subset C$ and
$P^x(T_C<\infty)=0$ for any $x\in E$. Let $m$ be a measure on
$(E,{\cal B})$. $D$ is called {\it $m$-essentially polar} if there
exists a set $C\in {\cal B}^n$ such that $D\subset C$ and
$P^m(T_C<\infty)=0$. Hereafter
$P^m(\cdot):=\int_EP^x(\cdot)m(dx)$.

Hunt's hypothesis (H) says that ``{\it every semipolar set of $X$ is
polar}". This hypothesis plays a crucial role in the potential
theory of (dual) Markov processes. It is known that
if $X$ is in duality with another standard process $\hat{X}$ on
$E$ with respect to a $\sigma$-finite reference measure $m$, then
(H) is equivalent to many potential principles for Markov
processes. Denote by $E^x$ the expectation with respect to $P^x$.  Let $\alpha>0$. A finite $\alpha$-excessive function $f$ on $E$ is called a regular potential provided that
$E^x\{e^{-\alpha T_n}f(X_{T_n})\}\rightarrow E^x\{e^{-\alpha
T}f(X_{T})\}$  for any $x\in E$ whenever $\{T_n\}$ is an increasing
sequence of stopping times with limit $T$. Denote by
$(U^{\alpha})_{\alpha>0}$ the resolvent operators for $X$.
\begin{itemize}
\item {\bf Bounded positivity principle $(P^*_{\alpha})$}: If $\nu$ is a finite signed measure
such that  $U^{\alpha}\nu$ is bounded, then $\nu U^{\alpha}\nu\geq 0$, where $\nu U^{\alpha}\nu:=\int_EU^{\alpha}\nu(x)\nu(dx)$.

\item {\bf Bounded energy principle $(E^*_{\alpha})$}: If $\nu$ is a finite measure with compact support
such that  $U^{\alpha}\nu$ is bounded, then $\nu$ does not charge
semipolar sets.

\item {\bf Bounded maximum principle $(M^*_{\alpha})$}: If $\nu$ is a finite measure
with compact support $K$ such that $U^{\alpha}\nu$ is bounded, then
$\sup\{U^{\alpha}\nu(x):x\in E\}=\sup\{U^{\alpha}\nu(x):x\in K\}$.

\item {\bf Bounded regularity principle $(R^*_{\alpha})$}: If $\nu$ is a finite measure
with compact support such that $U^{\alpha}\nu$ is bounded, then
$U^{\alpha}\nu$ is regular.

\end{itemize}

\begin{thm}\label{thm1} (Blumenthal and Getoor \cite{BG68,BG70}, Rao \cite{R77} and Fitzsimmons \cite{Fi90})
Assume that all 1-excessive (equivalently, all $\alpha$-excessive, $\alpha>0$) functions are lower
semicontinuous. Then
$$(P^*_{\alpha})\Leftrightarrow (E^*_{\alpha})
\Leftrightarrow(M^*_{\alpha})\Leftrightarrow(R^*_{\alpha})\Leftrightarrow ({\rm H}).
$$
\end{thm}

Hunt's hypothesis (H) is also equivalent to some other important properties of Markov processes. For example, (H) holds if and only if the fine and cofine topologies differ by
polar sets (Blumenthal and Getoor \cite[Proposition 4.1]{BG70} and Glover \cite[Theorem
2.2]{G83}); (H) holds if and only if every natural additive
functional of $X$ is a continuous additive functional
(Blumenthal and Getoor \cite[Chapter VI]{BG68}); (H) is equivalent to the dichotomy of capacity (Fitzsimmons and Kanda \cite{FK92}), which means that each compact set $K$ contains two disjoint sets with the same capacity as $K$.

In spite of its importance, (H) has been verified only in some special
situations. Some fifty years ago, Getoor conjectured that {\it
essentially all L\'{e}vy processes satisfy (H), except for some
extremely nonsymmetric cases like uniform motions}. This conjecture
stills remains open and is a major unsolved problem in the
potential theory for Markov processes.

The rest of this paper is organized as follows. In Section 2, we
survey the  existing results on Hunt's hypothesis (H) for Markov
processes and Getoor's conjecture for L\'{e}vy processes. In
Sections 3 and 4, we investigate (H) for multidimensional L\'{e}vy
processes from the viewpoints of projections and energy,
respectively. In Section 5, we present a few open questions for
further study.

\section{Survey on (H) for Markov processes}\setcounter{equation}{0}

In this section, we summarize the results that have  been obtained so
far for the validity of Hunt's Hypothesis  (H). We divide them into two parts: \S 2.1 (H) for L\'{e}vy processes and \S 2.2 (H) for Markov processes.

\subsection{(H) for L\'{e}vy processes}

Throughout this subsection, we let $(\Omega,{\cal F},P)$ be a probability space and $X=(X_t)_{t\ge 0}$ be an $\mathbf{R}^n$-valued L\'{e}vy process on $(\Omega,{\cal F},P)$ with L\'{e}vy-Khintchine exponent
$\psi$, i.e.,
\begin{eqnarray*}
E[\exp\{i\langle z,X_t\rangle\}]=\exp\{-t\psi(z)\},\ \ z\in
\mathbf{R}^n,t\ge 0.
\end{eqnarray*}
Hereafter we use  $E$ to denote the expectation {with respect to} $P$, and use $\langle\cdot,\cdot\rangle$ and $|\cdot|$ to denote the Euclidean inner product
and norm of $\mathbf{R}^n$, respectively. The classical L\'{e}vy-Khintchine formula tells us that
\begin{eqnarray*}
\psi(z)=i\langle a,z\rangle+\frac{1}{2}\langle z,Qz\rangle+\int_{\mathbf{R}^n}
\left(1-e^{i\langle z,x\rangle}+i\langle z,x\rangle 1_{\{|x|<1\}}\right)\mu(dx),
\end{eqnarray*}
where $a\in \mathbf{R}^n,Q$ is a symmetric nonnegative definite $n\times n$ matrix, and
$\mu$ is a measure (called the L\'evy measure) on $\mathbf{R}^n\backslash\{0\}$
satisfying $\int_{\mathbf{R}^n\backslash\{0\}} (1\wedge |x|^2)\mu(dx)<\infty$. For $x\in \mathbf{R}^n$,
we denote by $P^x$ the law of $x+X$ under $P$. In particular, $P^0=P$. Denote by $m_n$ the Lebesgue measure on $\mathbf{R}^n$.

We use Re$(\psi)$ and
Im$(\psi)$ to denote respectively the real and imaginary parts of $\psi$, and
use also $(a,Q,\mu)$ to denote $\psi$. Define
$$A:=1+{\rm Re}(\psi),\ \ B:=|1+\psi|.
$$
For a finite (positive) measure $\nu$ on $\mathbf{R}^n$, we denote
$$
\hat{\nu}(z):=\int_{\mathbf{R}^n}e^{i\langle z,x\rangle}\nu(dx).
$$
$\nu$ is said to have finite 1-energy if
$$
\int_{\mathbf{R}^n}\frac{A(z)}{B^2(z)}|\hat{\nu}(z)|^2dz<\infty.
$$
We use $\log$ to denote $\log_e$.

\subsubsection{Main results obtained before 1990}

Suppose that $X$ is a compound Poisson process. Then every $x\in \mathbf{R}^n$ is regular for $\{x\}$, i.e., $P^x(T_{\{x\}}=0)=1$. Hence only the empty set is a semipolar set and therefore (H) holds.

When $n=1$, Kesten \cite{Ke69} {(see also Bretagnolle \cite{Br71}) showed   that if $X$ is not a compound Poisson process, then  every $\{x\}$ is non-polar} if and only if
\begin{eqnarray}\label{Ke69-a}
\int_0^{\infty}\mbox{Re}([1+\psi(z)]^{-1})dz<\infty.
\end{eqnarray}
It follows that if condition (\ref{Ke69-a}) is fulfilled, then (H) holds is equivalent to that  only the empty set is a semipolar set.

Port and Stone \cite{PS69} proved that for the asymmetric
Cauchy process on the line every $x\in \mathbf{R}^n$ is regular for $\{x\}$. It follows that (H) holds in this case. Further, Blumenthal and Getoor
\cite{BG70} showed that all  stable processes with index
$\alpha\in (0,2)$ on the line satisfy (H). Kanda and Forst proved independently the following celebrated result.

\begin{thm}\label{thm-2.1-1} (Kanda \cite{Ka76} and Forst \cite{F75}) If $X$ has bounded continuous transition densities with respect to $m_n$ and  $|\mbox{Im} (\psi)|\leq
M(1+\mbox{Re}(\psi))$ for some positive constant $M$, then $X$ satisfies (H).
\end{thm}

Rao \cite{R77} gave a short proof of the above Kanda-Forst theorem under the weaker condition that $X$ has resolvent densities with respect to $m_n$. In fact, Rao proved that the  bounded maximum principle holds for $X$ in this case. By the Kanda-Forst theorem, we know that for $n\ge1$ all
stable processes with index $\alpha\neq 1$ satisfy (H). For the case $\alpha= 1$, under the additional assumption that the linear term vanishes, Kanda
\cite{Ka78} showed that (H) holds by virtue of the following result.

\begin{thm}\label{thm-2.1-2} (Kanda \cite[Theorem 1]{Ka78})
 Assume that $X$ has bounded continuous transition densities with respect to $m_n$. Then, a set contains
a semipolar subset for $X$ which is not polar for $X$, if and only if it contains a nonpolar compact subset $H$ for $X$ such that the $\alpha$-capacity of $H$ for $X$ is uniformly bounded as $\alpha\uparrow\infty$.
\end{thm}

By Hawkes \cite[Theorem 2.1]{Ha79}, we know that a L\'evy process has resolvent densities with respect to $m_n$ if and only if all 1-excessive functions are lower semicontinuous. In \cite{R88}, Rao gave a remarkable extension of the Kanda-Forst theorem.

\begin{thm}\label{thm-2.1-3} {\rm (Rao \cite[Theorem 2]{R88})} Assume that $X$ has resolvent densities with respect to $m_n$. Suppose there is an increasing function $f$ on $[1,\infty)$ such that
$ \int_N^{\infty}(\lambda f(\lambda))^{-1}d\lambda=\infty $ for
any $N\ge 1$ and  $|1+\psi|\le (1+{\rm Re}(\psi))f(1+{\rm
Re}(\psi))$. Then (H) holds.
\end{thm}

The function $f$ in Theorem \ref{thm-2.1-3} can be taken as: (i) $f\equiv M$ for some positive constant $M$ (now Rao's condition is reduced to the Kanda-Forst condition); (ii) $f(x)=\ln x$ for $x\geq 1$; (iii) $f(x)=\ln(1+\ln x)$ for $x\geq 1$.  By Theorem \ref{thm-2.1-3}, we find that the assumption  ``the linear term vanishes'' put in \cite{Ka78} can be removed. Hence all stable processes on $\mathbf{R}^n$ satisfy (H).

To prove Theorem \ref{thm-2.1-3}, Rao made use of the following result, which is a consequence of Rao \cite[Lemma 2.1 and Theorem 1.5]{R87}.

\begin{thm}\label{thm-2.1-4} {\rm (Rao \cite[Theorem 1]{R88})} Assume that $X$ has resolvent densities with respect to $m_n$.
Let $\nu$ be a finite measure of finite 1-energy.
Then
$$\lim_{\lambda\to\infty}\int_{\mathbf{R}^n}\frac{\lambda+{\rm Re}\psi(z)}{|\lambda+\psi(z)|^{2}}|\hat{\nu}(z)|^2dz
$$
exists. The limit is zero if and only if $U^1\mu$  is regular.
\end{thm}

\subsubsection{Main results obtained by our group}

 In this subsection, we introduce the main results obtained by our group on (H) for L\'{e}vy processes. We refer the reader to  Hu and Sun \cite{HS12}, Hu et al. \cite{HSZ15}, Hu and Sun \cite{HS16}, Hu and Sun \cite{HS18}, and Hu et al. \cite{HSW19} for more details.  The results are presented for three cases: (H) for general L\'{e}vy processes, (H) for one-dimensional L\'{e}vy processes, and (H) for subordinators.

\bigskip

\noindent{\bf (1) (H) for general L\'{e}vy processes}

\begin{thm}\label{thm-2.1-5} (\cite[Theorem 1.1]{HS12})
Suppose that $Q$ is non-degenerate, i.e., $Q$ is of full rank. Then:\\
(i) $X$ satisfies  (H);\\
(ii) The Kanda-Forst condition $|\mbox{Im} (\psi)|\leq
M(1+\mbox{Re}(\psi))$ holds for some positive constant $M$;\\
(iii) $X$ and $\tilde{X}$ have the same polar sets, where $\tilde{X}=X-\bar{X}$ with $\bar{X}$ being an independent copy of $X$.
\end{thm}

We first proved (ii) by showing that
\begin{equation}\label{add}
{\rm Re}\psi(z)\geq c\langle z,z\rangle,\ \  z\in \mathbf{R}^n,
\end{equation}
 for some positive constant $c$ and using the inequality that $|t-\sin t|\le t^2/2$, $t\in \mathbf{R}$. By (\ref{add}) and Hartman and Wintner \cite{HW42}, one finds that $X$ has bounded continuous transition densities. Hence (i) holds by the Kanda-Forst theorem. (iii) is a consequence of (ii), Kanda \cite[Theorem 1]{Ka76} (or Hawkes \cite[Theorems 2.1 and 3.3]{Ha79}), and the fact that $X$ has bounded continuous transition densities.

Denote $b:=-a$ and $\mu_1:=\mu|_{\mathbf{R}^n\backslash
\sqrt{Q}\mathbf{R}^n}$. If
$\int_{\{|x|<1\}}|x|\mu_1(dx)<\infty$, we set $b':=b-\int_{\{|x|<1\}}x\mu_1(dx)$. Define the following solution condition:

{\it  (S)\ \ The equation $\sqrt{Q}y=b',\ y\in \mathbf{R}^n$,
has at least one solution.}

\noindent Note that (S) means that $b'\in \mathcal{R}(\sqrt{Q}):=\{\sqrt{Q}y: y\in \mathbf{R}^n\}$.

\begin{thm}\label{thm-2.1-6} (\cite[Theorem 1.2]{HS12})
Suppose that $\mu({\mathbf{R}^n\backslash
\sqrt{Q}\mathbf{R}^n})<\infty$. Then, the following three claims are
equivalent:

(i) $X$ satisfies  (H);

(ii)  (S) holds;

(iii) The Kanda-Forst condition $|{\rm Im} (\psi)|\leq M(1+{\rm Re}(\psi))$ holds for some
positive constant $M$.

\end{thm}

The L\'{e}vy-It\^{o} decomposition, orthogonal transformation and strong Markov property of L\'{e}vy processes, and some properties of compound Poisson processes have been used to prove Theorem \ref{thm-2.1-6}.

\begin{pro}\label{pro-2.1-7} (\cite[Proposition 1.5]{HS12})
Suppose that $X$ has bounded continuous transition densities,
 and $X$ and $\tilde{X}$ have the same polar sets, where $\tilde{X}$ is defined in Theorem \ref{thm-2.1-5}(iii). Then $X$ satisfies  (H).
\end{pro}

To prove the above proposition, we used mainly an idea given in the proof of Kanda \cite[Theorem 2]{Ka78} and the comparison inequality for capacities given in Kanda \cite{Ka76} (cf. also Hawkes \cite{Ha79}).

In \cite{R88}, Rao mentioned that his condition $|1+\psi|\le (1+{\rm Re}(\psi))f(1+{\rm Re}(\psi))$, equivalently, $|{\rm Im}(\psi)|\le (1+{\rm Re}(\psi))f(1+{\rm Re}(\psi))$, is not far from being necessary for the validity of (H). The following result tells us that Rao's condition can be relaxed.

\begin{thm}\label{thm-2.1-8} (\cite[Theorem 4.5]{HSZ15}) Assume that all 1-excessive functions are lower semicontinuous. Then
(H) holds if the following extended Kanda-Forst-Rao condition holds:

(EKFR) There are two measurable functions $\psi_1$ and $\psi_2$ on $\mathbf{R}^n$ such that $\rm{Im}(\psi)=\psi_1+\psi_2$, $|\psi_1|\leq Af(A)$, and
\begin{eqnarray*}
\int_{{\mathbf{R}^n}}\frac{|\psi_2(z)|}{(1+{\rm Re}\psi(z))^2+({\rm Im}\psi(z))^2}dz<\infty,
\end{eqnarray*}
where  $f$ is a positive increasing function on $[1,\infty)$ such that $
\int_N^{\infty}(\lambda f(\lambda))^{-1}d\lambda=\infty $ for some
$N\ge 1$.
\end{thm}

By virtue of Theorem \ref{thm-2.1-8}, we constructed a class of one-dimensional L\'{e}vy processes satisfying (H) in \cite[Example 4.8]{HSZ15}. These L\'{e}vy processes have sufficient number of small jumps and no restriction is put on $a$ or $Q$. As another application of Theorem \ref{thm-2.1-8}, we have the following result.

\begin{pro}\label{pro-2.1-23}(\cite[Proposition 4.10]{HSZ15}) Let  $X$ be a L\'{e}vy process on $\mathbf{R}$ such that  all 1-excessive functions are lower semicontinuous. Suppose that $$
\liminf_{|z|\to \infty}\frac{|\psi(z)|}{|z|\log^{1+\gamma}|z|}>0
$$
for some constant $\gamma>0$. Then (H) holds.
\end{pro}

To prove Theorem \ref{thm-2.1-8}, we used the following necessary and sufficient condition for (H).

\begin{thm}\label{thm-2.1-9} (\cite[Theorem 4.3]{HSZ15}) Assume that all 1-excessive functions are lower semicontinuous. Let $f$ be a positive increasing
function on $[1,\infty)$ such that $ \int_N^{\infty}(\lambda
f(\lambda))^{-1}d\lambda=\infty $ for some $N\ge 1$. Then (H)
holds if and only if
\begin{eqnarray*}
\lim_{\lambda\rightarrow\infty}\sum_{k=1}^{\infty}\int_{\left\{B(z)>A(z)f(A(z)),\,k\le \frac{|{\rm Im}\psi(z)|}{A(z)}<k+1,\,A(z)\le\lambda<(k+1)|{\rm Im}\psi(z)|\right\}}\frac{\lambda}{\lambda^2+({\rm Im}\psi(z))^2}|\hat \nu(z)|^2dz=0
\end{eqnarray*}
for any finite measure $\nu$ with compact support such that $U^1\nu$ is bounded.
\end{thm}

The proof of Theorem \ref{thm-2.1-9} is based on a key lemma (\cite[Lemma 4.2]{HSZ15}), which is obtained by using \cite[Theorems 1 and 2]{R88}.

The following two theorems provide new necessary and sufficient conditions for the validity of
(H) for L\'{e}vy processes. Different from the classical Kanda-Forst condition
and Rao's condition, our conditions only require that ${\rm Im}(\psi)$ is partially
well-controlled by $1+{\rm Re}(\psi)$. The weaker conditions are fulfilled by more general
L\'{e}vy processes and reveal the more essential reason for the validity of (H) (see \cite[Section 3]{HS16} for examples).

\begin{thm}\label{thm-2.1-10} (\cite[Theorem 2.3]{HS16}) Assume that  all 1-excessive functions are lower semicontinuous.

(i) $X$ satisfies (H)  if the following condition holds:
\begin{eqnarray}\label{zeta}
& &{ Condition\ (C^{\log})}:\ { For\ any\ finite\
measure}\ \nu\ { on}\ \mathbf{R}^n\ { of\ finite\ 1{\textrm{-}}energy,\ there}\nonumber\\
& &{exist\ a\ constant}\ \varsigma>1\  and\ a\ sequence\
\ \{y_k\uparrow\infty\}\ { such\ that}\ y_1>1\ and\nonumber\\
& &\ \ \ \ \ \ \ \ \ \ \ \ \ \ \sum_{k=1}^{\infty}\int_{\{y_k\le
B(z)< (y_k)^{\varsigma}\}} \frac{1}{B(z)\log(B(z))}|\hat
\nu(z)|^2dz<\infty.
\end{eqnarray}

(ii) Suppose $X$ satisfies (H). Then, for\ any\ finite\ measure
$\nu$ on $\mathbf{R}^n$  of\ finite\ 1-energy and any
$\varsigma>1$, there exists a sequence\ $\{y_k\uparrow\infty\}$
such\ that $y_1>1$ and (\ref{zeta})\ {holds}.
\end{thm}

Let $\varsigma>1$ be a constant. We define
$N^{\varsigma}_x:=\varsigma^{(\varsigma^x)}\ \ {\rm for}\
x\ge 0.$

\begin{thm}\label{thm-2.1-11} (\cite[Theorem 2.4]{HS16}) Assume that all 1-excessive functions are lower semicontinuous.

 (i)
$X$ satisfies (H)  if the following condition holds:
\begin{eqnarray}\label{new345}
& &{ Condition\ (C^{\log\log})}:\ { For\ any\ finite\
measure}\ \nu\ { on}\ \mathbf{R}^n\ { of\ finite\ 1{\textrm{-}}energy,}\nonumber\\
& &\ \ \ \ \ {there\ exist\ a\ constant}\ \varsigma>1\ { and\ a\
sequence\ of\ positive\ numbers}\
\{x_k\}\nonumber\\
& &\ \ \ \ \ { such\ that}\ N^{\varsigma}_{x_1}>e,\
x_k+1<x_{k+1},\ k\in
\mathbf{N},\ \sum_{k=1}^{\infty}\frac{1}{x_k}=\infty,\ { and}\nonumber\\
& &\ \ \ \ \ \ \sum_{k=1}^{\infty}\int_{\{N^\varsigma_{x_k}\le
B(z)< N^\varsigma_{x_k+1}\}} \frac{1}{B(z)\log
(B(z))[\log\log(B(z))]}|\hat \nu(z)|^2dz<\infty.
\end{eqnarray}

(ii) Suppose $X$ satisfies (H). Then, for any finite measure $\nu$
on $\mathbf{R}^n$ of finite 1-energy and any $\varsigma>1$, there
exists a sequence of positive numbers $\{x_k\}$ such that $
N^{\varsigma}_{x_1}>e$, $x_k+1<x_{k+1}$, $k\in \mathbf{N}$, $
\sum_{k=1}^{\infty}\frac{1}{x_k}=\infty$, and (\ref{new345})\
{holds}.
\end{thm}

The proofs of the above two theorems rely on the following characterization for (H).

\begin{pro} \label{pro-2.1-12} (\cite[Proposition 2.2]{HS16}) Assume that all 1-excessive functions are lower semicontinuous. Then (H) holds if and only if
\begin{eqnarray*}
\lim_{\lambda\to\infty}\int_{\mathbf{R}^n}\frac{\lambda}{\lambda^2+B^2(z)}|\hat\nu(z)|^2dz=0
\end{eqnarray*}
for any finite measure $\nu$ of finite 1-energy.
\end{pro}

Motivated by exploration of Getoor's conjecture for one-dimensional L\'{e}vy processes (see \cite[Section 2.1]{HS18}),  we considered in \cite{HS18} Hunt's hypothesis (H) for the sum of two independent L\'{e}vy processes.

\begin{thm}\label{thm-2.1-13} (\cite[Theorem 3.1]{HS18})
Let $X_1$ and $X_2$ be two independent L\'evy processes on
$\mathbf{R}^n$. If $X_1$ satisfies (H) and $X_2$ is a compound
Poisson process, then $X_1+X_2$ satisfies (H).
\end{thm}

Hereafter we say that a L\'{e}vy process with L\'{e}vy-Khintcine exponent $(a,Q,\mu)$ satisfies {\it condition (S)} if
 $\mu({\mathbf{R}^n\backslash \sqrt{Q}\mathbf{R}^n})<\infty$
and  the equation $ \sqrt{Q}y=-a-\int_{\{\mathbf{R}^n\backslash
\sqrt{Q}\mathbf{R}^n\}}x1_{\{|x|<1\}}\mu(dx) $ has at least one
solution $y\in \mathbf{R}^n$.
\begin{thm}\label{thm-2.1-14}(\cite[Theorem 3.2]{HS18})
Let $X_1$ and $X_2$ be  two independent L\'evy processes on
$\mathbf{R}^n$. If both $X_1$ and $X_2$ satisfy condition (S),
then $X_1+X_2$ satisfies (H).
\end{thm}

To show Theorem \ref{thm-2.1-13}, we considered projections for L\'{e}vy processes (see \cite[Lemma 3.4]{HS18}) and used an idea in the proof of \cite[Theorem 1.2]{HS12} (see \cite[Lemma 3.6]{HS18}). To show Theorem \ref{thm-2.1-14}, we proved a lemma for general symmetric nonnegative matrices (see \cite[Lemma 3.7]{HS18}).

If L\'{e}vy processes have resolvent densities, we have the following result on the validity of (H).

\begin{thm}\label{thm-2.1-15} (\cite[Theorem 4.1]{HS18})
Assume that $X_1$ and $X_2$ are two
independent L\'{e}vy processes on $\mathbf{R}^n$ such that
$X_1+X_2$ has resolvent densities with respect to $m_n$. Denote by $\psi_1$ and $\psi_2$ the L\'{e}vy-Khintchine exponents
of $X_1$ and $X_2$, respectively.
Suppose that

(i) $X_1$ has resolvent densities with respect to $m_n$ and
satisfies (H).

(ii) Any finite measure $\nu$ of finite 1-energy with respect to $X_1+X_2$
has finite 1-energy with respect to $X_1$.

(iii) There exists a constant $c>0$ such that $|{\rm Im}(\psi_2)|\leq c(1+{\rm Re}(\psi_1)+{\rm Re}(\psi_2)).$\\
Then $X_1+X_2$ satisfies (H).
\end{thm}
For condition (ii) of Theorem \ref{thm-2.1-15}, we refer the reader to \cite[Proposition 4.2]{HS18} for some sufficient conditions.

Before ending this subsection, we present a result which implies that big jumps have no effect on the validity of (H) for any
L\'{e}vy process.

\begin{thm}\label{thm-2.1-16}(\cite[Proposition 4.11]{HSW19})
 Suppose  that $\mu_1$
is a finite measure on $\mathbf{R}^n\backslash\{0\}$ such that
$\mu_1\leq \mu$. Denote $\mu':=\mu-\mu_1$ and let $X'$ be a
L\'{e}vy process on $\mathbf{R}^n$ with L\'{e}vy-Khintchine
exponent $(a',Q,\mu')$, where $a':=a+\int_{\{|x|< 1\}}x\mu_1(dx).$
Then,

(i) $X$ and $X'$ have same semipolar sets.

(ii) $X$ and $X'$ have same $m_n$-essentially polar sets.

(iii) $X$ satisfies (H) if and only if $X'$ satisfies (H).

(iv) $X$ satisfies $(H_{m_n})$ if and only if $X'$ satisfies $(H_{m_n})$, where $(H_{m_n})$ means that every semipolar set is $m_n$-essentially polar.
\end{thm}

We proved Theorem \ref{thm-2.1-16} (i) and (ii) in \cite[Theorem 2.1]{HSZ15}. (iv) is a direct consequence of (i) and (ii), and (iii) is based on \cite[Theorem 3.1]{HS18} and \cite[Theorem 1.1, Corollary 4.10]{HSW19}.

\bigskip

\noindent{\bf (2) (H) for one-dimensional L\'{e}vy processes}

\noindent In this part, we assume that $X$ is a one-dimensional L\'{e}vy process with L\'{e}vy-Khintchine exponent $(a,Q,\mu)$. Denote by $\mu_{+}$ and $\mu_{-}$ the restriction of $\mu$ on $(0,\infty)$ and $(-\infty,0)$, respectively.
 Let
$\bar{\mu}_{-}$ be the image measure of $\mu_{-}$ under the map
$
x\mapsto -x, \ \forall x\in (-\infty,0).
$
The following result extends  Kesten \cite[Theorem 1(f)]{Ke69}.

\begin{thm}\label{thm-2.1-17} (\cite[Theorem 2.2]{HS18})
Suppose that $Q=0$ and  $\int_{0}^{\infty}(1\wedge x
)\mu_+(dx)=\infty$. If there exist $\delta\in (0,1), k\in [0,1)$,
and a measure $\nu$ on $\mathbf{R}^+$ satisfying
$\int_0^{\delta}x\nu(dx)<\infty$, such that
$\bar{\mu}_{-}\leq k\mu_++\nu.$ Then $X$ satisfies (H).
\end{thm}

 The basic idea of the proof for Theorem  \ref{thm-2.1-17} is to use Kesten's criterion  (\ref{Ke69-a}) and Bretagnolle's beautiful characterization of one-dimensional L\'{e}vy processes (see \cite[Theorem 8]{Br71}).

We now give a novel condition on the L\'{e}vy measure $\mu$ which implies (H) for a large class of
one-dimensional L\'{e}vy processes.

\begin{thm}\label{thm-2.1-18}(\cite[Theorem 2.2]{HS18})
If
\begin{eqnarray}\label{thm-2.1-18-a}
\liminf_{\varepsilon\rightarrow
0}\frac{\int_{-\varepsilon}^{\varepsilon}x^2\mu(dx)}{\varepsilon/|\log
\varepsilon|}>0,
\end{eqnarray}
then $X$ satisfies (H).
\end{thm}

The proof of Theorem \ref{thm-2.1-18} is based on the following characterization for (H).

\begin{pro}\label{pro-2.1-19} (\cite[Proposition 2.4]{HS18}) Suppose that $X$ is a L\'{e}vy process on $\mathbf{R}^n$  which has resolvent densities with respect to $m_n$. Let
 $f$ be a positive increasing
function on $[1,\infty)$ such that $ \int_N^{\infty}(\lambda
f(\lambda))^{-1}d\lambda=\infty $ for some $N\ge 1$. Then (H)
holds for $X$ if and only if
\begin{eqnarray*}
\lim_{\lambda\rightarrow\infty}\int_{\{B(z)>A(z)f(A(z))\}}\frac{\lambda}{\lambda^2+B^2(z)}|\hat
\nu(z)|^2dz=0
\end{eqnarray*}
for any finite measure $\nu$ of finite 1-energy.
\end{pro}

\begin{rem}\label{rem-2.1-20} Note that, different from most existing sufficient conditions for (H), our condition (\ref{thm-2.1-18-a}) does not require any controllability of ${\rm Im}(\psi)$ by $1+{\rm Re}(\psi)$. Define the measure $\xi$ on $\mathbf{R}$ by
$$
\xi(dx) := (1\wedge |x|)\mu(dx),\ \ x\in \mathbf{R}.
$$
Condition (\ref{thm-2.1-18-a}) is
slightly stronger than $\xi$ is an infinite
measure on $\mathbf{R}$. We refer the reader to \cite[Remark  2.5]{HS18} for more details.
\end{rem}

From the proof of Theorem \ref{thm-2.1-18}, we we can see that the following result extending \cite[Theorem 4.7]{HSZ15} holds.

\begin{pro}\label{pro-2.1-21} (\cite[Proposition 2.6]{HS18})
If
$$
\liminf_{|z|\to\infty}\frac{{\rm Re}\psi(z)}{|z|/\log|z|}>0,
$$
then $X$ satisfies (H).
\end{pro}

Following the proof of  Theorem \ref{thm-2.1-18}, we can prove the following result.

\begin{pro}\label{pro-2.1-22} (\cite[Proposition 2.7]{HS18})
If
$$
\liminf_{\varepsilon\to
0}\frac{\int_{-\varepsilon}^{\varepsilon}x^2\mu(dx)}{\frac{\varepsilon}{|\log
\varepsilon|[\log|\log\varepsilon|]}}>0,
$$
 then $X$ satisfies (H).
\end{pro}

\noindent {\bf (3) (H) for subordinators}

\noindent $X$ is called a subordinator if it is a one-dimensional increasing L\'{e}vy process. Subordinators are a very important class of L\'{e}vy processes. Let $X$ be a subordinator. Then its L\'{e}vy-Khintchine exponent $\psi$ can be expressed by
$$
\psi(z)=-idz+\int_0^{\infty}{\left(1-e^{izx}\right)}\mu(dx),\
z\in \mathbf{R},
$$
where $d\geq 0$ (called the drift coefficient) and $\mu$ satisfies
$\int_0^{\infty}(1\wedge x)\mu(dx)<\infty$.

\begin{pro}\label{pro-2.1-24} (\cite[Proposition 1.6]{HS12})
If $X$ is a subordinator and satisfies (H), then $d=0$.
\end{pro}

\noindent Proposition \ref{pro-2.1-24} can be extended to the high-dimensional case, see Proposition 3.2 below.

A natural question is: {\it if $X$ is a pure jump subordinator, i.e., $d=0$, must $X$ satisfy (H)?}  Up to now, it is still unknown if the answer is yes or no. In the following, we first show that some particular subordinators satisfy (H).

Recall that the potential measure $U$ of $X$ is defined by
$$
U(A)=E\left[\int_0^{\infty}1_A(X_t)dt\right],\ \  A\subset [0,\infty).
$$
$X$ is called a {\it  special subordinator} if $U|_{(0,\infty)}$ has a decreasing density with respect to the Lebesgue measure.

\begin{thm}\label{thm-2.1-25} (\cite[Theorem 3.3]{HSZ15})
Let $X$ be a special subordinator.  Then $X$ satisfies (H) if and
only if $d=0$.
\end{thm}

\begin{defi}\label{defi-2.1-26} (\cite[Definition 3.4]{HSZ15})
Let $X$ be a subordinator with drift $0$ and L\'{e}vy measure
$\mu$. We call $X$ a {locally quasi-stable subordinator} if there
exist a stable subordinator $S$ with L\'{e}vy measure $\mu_S$,
positive constants $c_1,c_2,\delta$, and finite measures $\mu_1$
and $\mu_2$ on $(0,\delta)$ such that
$$c_1\mu_S-\mu_1\le
\mu\le c_2\mu_S+\mu_2\ \ {\rm on}\ (0,\delta). $$
\end{defi}

\begin{pro}\label{pro-2.1-27} (\cite[Proposition 3.5]{HSZ15})
Any locally quasi-stable subordinator satisfies (H).
\end{pro}

We refer the reader to \cite[Section 3.3]{HSZ15} and \cite[Example 4.10]{HS18} for more examples on subordinators satisfying (H).

In \cite[Section 5]{HSZ15}, we constructed a type of subordinators
that does not satisfy Rao's condition. So far we have not been able to prove or disprove that (H) holds for the subordinators. The example suggests that maybe completely new ideas and methods are needed for resolving Getoor's conjecture.

\begin{defi}\label{add2} (\cite[Definition 4.1]{HS16})
Let $0 <\alpha<\beta< 1$. A pure jump subordinator $X$ is said to be of
type-$(\alpha,\beta)$ if the L\'{e}vy measure of $X$ has density, which is denoted by $\rho$, and there
exists a constant $c > 1$ such that
$$
\frac{1}{cx^{1+\alpha}}\le \rho(x)\le \frac{c}{x^{1+\beta}},\ \ \forall x\in (0, 1].
$$
\end{defi}

Up to now it is still unknown if any pure jump subordinator of type-$(\alpha,\beta)$ satisfies (H). But we have proved the following result based on Theorem \ref{thm-2.1-10}.

\begin{thm} (\cite[Theorem 4.2]{HS16})
Any pure jump subordinator of type-$(\alpha,\beta)$ can be decomposed into
the summation of two independent pure jump subordinators of type-$(\alpha,\beta)$ such that
both of them satisfy (H).
\end{thm}

\subsection{(H) for Markov processes}

In this subsection, we assume that  $E$ is a locally compact space with a countable base and $X=(X_t, P^x)$ is a standard Markov process on $E$.

Suppose that $X$ is   associated with a (not necessarily symmetric) regular Dirichlet form on $L^2(E;m)$, where $m$ is a Radon measure on $E$. Silverstein \cite{Si77} proved that any semipolar set for $X$ is $m$-essentially polar. This result plays a very important role in the theory of Dirichlet forms. For example, it is used to prove the relationship between orthogonal projections and hitting distributions (cf. \cite[Theorem 4.3.1]{FOT} and its proof).  Fitzsimmons \cite{Fi01} extended the result to the
semi-Dirichlet forms setting and  Han et al. \cite{HMS11} extended it to the
positivity-preserving forms setting.

In \cite{GR88}, Glover and Rao gave a sufficient condition  for nonsymmetric Hunt processes to satisfy (H).  In
\cite{Fi14}, Fitzsimmons  showed that Gross's Bwownian motion,
which is an infinite-dimensional L\'evy process, fails to satisfy
(H). In \cite{HN16}, Hansen and Netuka showed that (H) holds if
there exists a Green function $G>0$ which locally satisfies the
triangle inequality $G(x,z)\wedge G(y,z)\leq CG(x,y)$, where $C$ is a positive constant.

In \cite{HSW19}, we investigated the invariance of  (H) for Markov processes under two classes of transformations, which are change of measure and
subordination. Before stating our results, we give some notation. We fix an isolated point
$\Delta$ which is not in $E$ and write $E_{\Delta}=E\cup
\{\Delta\}$. Consider the following objects:

(i) $\Omega$ is a set and $\omega_{\Delta}$ is a distinguished point of $\Omega$.

(ii) For $0\le t\le\infty$, $Z_t:\Omega \rightarrow E_{\Delta}$ is
a map such that if $Z_t(\omega)=\Delta$  then $Z_s(\omega)=\Delta$
for all $s\ge t$, $Z_{\infty}(\omega)=\Delta$ for all
$\omega\in\Omega$, and $Z_0(\omega_{\Delta})=\Delta$.

(iii) For $0\le t\le\infty$, $\theta_t:\Omega\rightarrow\Omega$ is a map such that $Z_s\circ \theta_t=Z_{s+t}$ for all $s,t\in [0,\infty]$, and $\theta_{\infty}\omega=\omega_{\Delta}$ for all $\omega\in\Omega$.

We define in $\Omega$ the $\sigma$-algebras
$\mathcal{F}^0=\sigma(Z_t: t\in [0,\infty])$ and
$\mathcal{F}^0_t=\sigma(Z_s: s\leq t)$ for $0\le t<\infty$. Denote
$$
\zeta(\omega)=\inf\{\omega:Z_t(\omega)=\Delta\},\ \ \omega\in \Omega.
$$
Let $m$ be a measure on $(E,{\cal B})$. We define
$$(H_m):\ \mbox{every semipolar set is}\ m\mbox{-essentially\ polar}.$$
Note that if a standard process $X$ has resolvent densities with
respect to $m$, then $X$ satisfies (H) if and only if $X$
satisfies $(H_m)$ (cf. \cite[Propositions II.2.8 and
II.3.2]{BG68}).

\begin{thm}\label{thm-2.2-2} (\cite[Theorem 1.1]{HSW19})
Let $X=(\Omega,\mathcal{M}^X,\mathcal{M}^X_t,Z_t,\theta_t,P^x)$ and
$Y=(\Omega,\mathcal{M}^Y,\mathcal{M}^Y_t,Z_t,$ $\theta_t,Q^x)$ be two
standard processes on $E$ such that $\mathcal{M}^X\cap \mathcal{M}^Y\supset \mathcal{F}^0$ and $\mathcal{M}^X_t\cap \mathcal{M}^Y_t\supset \mathcal{F}_t^0$ for $0\le t<\infty$.

(i) Suppose that $X$ satisfies (H) and for any $x\in E$ and $t>0$, $Q^x|_{\cal{F}_t^0}$ is absolutely
continuous with respect to $P^x|_{\cal{F}_t^0}$ on $\{t<\zeta\}$. Then $Y$ satisfies (H).

(ii) Suppose that $X$ satisfies ($H_m$) for some measure $m$ on $(E,{\cal B})$ and for any $x\in E$ and $t>0$, $Q^x|_{\cal{F}_t^0}$ is absolutely
continuous with respect to $P^x|_{\cal{F}_t^0}$ on $\{t<\zeta\}$. Then $Y$ satisfies ($H_m$).
\end{thm}

Let $X=(X_t)_{t\ge 0}$ be a standard
process on $E$ and $\tau$ be a subordinator  which is independent
of $X$. The standard process $(X_{\tau_t})_{t\ge 0}$ is called the
subordinated process of $(X_t)_{t\ge 0}$. The idea of subordination
originated from Bochner (cf. \cite{Bo55}). Our next result is
motivated by the following remarkable theorem of Glover and Rao.
\begin{thm} \label{thm-3}  (Glover and Rao \cite{GR86}) Let $(X_t)_{t\ge 0}$
be a standard process on $E$ and $(\tau_t)_{t\ge 0}$ be a subordinator
which is independent of $X$ and satisfies (H). Then $(X_{\tau_t})_{t\ge 0}$
satisfies (H).
\end{thm}

Now we present our result on the equivalence between (H) for $X$
and (H) for its time changed process.

\begin{thm}\label{thm-2.2-3} (\cite[Theorem 1.3]{HSW19})
Let $(X_t)_{t\ge 0}$ be a standard process on $E$ and $m$ be a measure on
$(E,{\cal B})$. Then,

(i) $(X_t)_{t\ge 0}$ satisfies (H) if and only if $(X_{\tau_t})_{t\ge 0}$ satisfies (H)
for some (and hence any) subordinator $(\tau_t)_{t\ge 0}$ which is
independent of $(X_t)_{t\ge 0}$ and has a positive drift coefficient.

(ii) $(X_t)_{t\ge 0}$ satisfies ($H_m$) if and only if $(X_{\tau_t})_{t\ge 0}$ satisfies ($H_m$)
for some (and hence any) subordinator $(\tau_t)_{t\ge 0}$ which is
independent of $(X_t)_{t\ge 0}$ and has a positive drift coefficient.
\end{thm}

The proof of Theorem \ref{thm-2.2-2} is based on two lemmas (\cite[Lemmas 2.1 and 2.2]{HSW19}) and Blumenthal's 0-1 law. The proof of Theorem \ref{thm-2.2-3} is based on \cite[Lemma 2.1]{HSW19} and Bertoin \cite[Theorem III.5]{B96}. We refer the reader to \cite[Theorems 4.3, 4.7, 4.13 and Proposition 5.1]{HSW19} for applications of Theorems \ref{thm-2.2-2} and \ref{thm-2.2-3}.

\section{(H) for multidimensional L\'{e}vy processes: projections}\setcounter{equation}{0}

In this section, we investigate  (H) for multidimensional L\'{e}vy processes from the  viewpoint of projections. Throughout this section, we assume that $n>1$ except in Proposition 3.2 below and $X=(X_t)_{t\ge 0}$ is a L\'{e}vy process on  $\mathbf{R}^n$ with L\'{e}vy-Khintchine exponent $(a,Q,\mu)$. For a subspace $A$  of $\mathbf{R}^n$, we use $A^{\perp}$ to denote its orthogonal complement space.

\subsection{A lemma on projections and applications}\setcounter{equation}{0}

In \cite{HS18}, we  proved the following result.

\begin{lem}\label{lem-3.1} (\cite[Lemma 3.4]{HS18})
Suppose that $X$ satisfies  (H). Then for any nonempty proper subspace $A$ of $\mathbf{R}^n$, the projection process $Y=(Y_t)_{t\ge 0}$ of $X$ on $A$  satisfies (H).
\end{lem}

As an application of Lemma \ref{lem-3.1}, we proved in \cite{HS18} Theorem \ref{thm-2.1-13} of  Section 2. As another application of Lemma \ref{lem-3.1}, we proved in \cite{HSW19} the following result, which extends Proposition \ref{pro-2.1-24}.

\begin{pro}\label{pro-3.2} (\cite[Proposition 5.3]{HSW19})
 Let $X$ be a L\'{e}vy process on  $\mathbf{R}^n\,(n\geq 1)$ with L\'{e}vy-Khintchine exponent $(a,0,\mu)$ satisfying $\int_{\mathbf{R}^n} (|x|\wedge 1)\mu(dx)<\infty$. If $X$ satisfies (H), then its drift coefficient equals zero.
\end{pro}

Proposition \ref{pro-3.2} can be further extended as follows.

\begin{pro}\label{pro-3.3}
Suppose that $Q$ is degenerate.  Let $Y=(Y_t)_{t\ge 0}$ be the projection process of $X$ on $(\sqrt{Q}\mathbf{R}^n)^{\perp}$. Denote by $P_2$   the projection operator  from $\mathbf{R}^n$ to $(\sqrt{Q}\mathbf{R}^n)^{\perp}$ and denote by $\mu_{P_2}$  the image measure of $\mu$ under $P_2$. Assume that
$$
\int_{(\sqrt{Q}\mathbf{R}^n)^{\perp}}( 1\wedge |x|)\mu_{P_2}(dx)<\infty.
$$
 If $X$ satisfies (H), then the drift coefficient of $Y$ equals zero.
\end{pro}

Proposition \ref{pro-3.3} is a direct consequence of Lemma  \ref{lem-3.1}, Proposition \ref{pro-3.2} and the following lemma.

\begin{lem}\label{lem-3.4}
Suppose that $Q$ is degenerate.  Let $Y=(Y_t)_{t\ge 0}$ be the
projection process of $X$ on $(\sqrt{Q}\mathbf{R}^n)^{\perp}$.
Denote by $P_2$   the projection operator  from $\mathbf{R}^n$ to
$(\sqrt{Q}\mathbf{R}^n)^{\perp}$ and denote by $\mu_{P_2}$  the
image measure of $\mu$ under $P_2$. Define
$$
a'=P_2a+\int_{\mathbf{R}^n}P_2x\left(1_{\{|x|<1\}}-1_{\{|P_2x|<1\}}\right)\mu(dx).
$$
Then  the L\'{e}vy-Khintchine exponent
of $Y$ is $(a',0,\mu_{P_2})$.
\end{lem}
{\bf Proof.} We use $k$ to denote the rank of $Q$ and assume without loss of generality that $1\leq k<n$. Then  there exists an orthogonal matrix
$O$ such that
\begin{eqnarray}\label{aa}
OQO^T=diag(\lambda_1,\dots,\lambda_n):=D,
\end{eqnarray}
where $\lambda_1\geq \lambda_2\geq \cdots\geq \lambda_k> 0$, $\lambda_i=0$ for $i=k+1,\dots,n$, and $O^T$ denotes the transpose of $O$. Define $U=O^Tdiag(\sqrt{\lambda_1},\dots,\sqrt{\lambda_n})O$. Then $Q=U^2$ and hence $\sqrt{Q}=U$.

Let $E_k=diag(1,\dots,1,0,\dots,0)$ be the diagonal matrix with the first $k$ elements being 1. Then, we have that
\begin{equation}\label{Q}
\sqrt{Q}\mathbf{R}^n=U\mathbf{R}^n=O^Tdiag(\sqrt{\lambda_1},\dots,\sqrt{\lambda_n})O\mathbf{R}^n=O^TE_kO\mathbf{R}^n.
\end{equation}
Define
$$
P_1=O^TE_kO.
$$
Then (\ref{Q}) implies  that $P_1$ is  the projection operator  from $\mathbf{R}^n$ to $\sqrt{Q}\mathbf{R}^n$. Define
\begin{equation}\label{bb}
P_2=I_n-P_1=O^T(I_n-E_k)O,
\end{equation}
where $I_n$ is the identity operator on $\mathbf{R}^n$. Then, $P_2$ is  the projection operator  from $\mathbf{R}^n$ to $(\sqrt{Q}\mathbf{R}^n)^{\perp}$.

Now we compute the L\'{e}vy-Khintchine exponent
of $Y$. Note that $Y_t=P_2X_t$ for $t\ge 0$. For $z\in \mathbf{R}^n$, we have
\begin{eqnarray*}
&&E\left[e^{i\langle z,P_2X_1\rangle}\right]=E\left[e^{i\langle P_2z,X_1\rangle}\right]\\
&&=\exp\left[-\left(i\langle a,P_2z\rangle+\frac{1}{2}\langle P_2z,QP_2z\rangle+\int_{\mathbf{R}^n}\left(1-e^{i\langle P_2z,x\rangle}+i\langle P_2z,x\rangle1_{\{|x|<1\}}\right)\mu(dx)\right)\right]\\
&&=\exp\left[-\left(i\langle P_2a,z\rangle+\frac{1}{2}\langle z,P_2QP_2z\rangle+\int_{\mathbf{R}^n}\left(1-e^{i\langle z,P_2x\rangle}+i\langle z,P_2x\rangle1_{\{|P_2x|<1\}}\right.\right.\right.\\
&&\quad\quad\quad\quad\left.\left.\left.+i\langle z,P_2x\rangle\left(1_{\{|x|<1\}}-1_{\{|P_2x|<1\}}\right)\right)\mu(dx)\right)\right]\\
&&=\exp\left[-\left(i\langle a',z\rangle+\frac{1}{2}\langle z,P_2QP_2z\rangle+\int_{\mathbf{R}^n}\left(1-e^{i\langle z,y\rangle}+i\langle z,y\rangle1_{\{|y|<1\}}\right)\mu_{P_2}(dy)\right)\right].
\end{eqnarray*}
By (\ref{aa}) and (\ref{bb}), we get
\begin{eqnarray*}
P_2QP_2&=&O^T(I_n-E_k)OO^TDOO^T(I_n-E_k)O=0.
\end{eqnarray*}
Therefore, the L\'{e}vy-Khintchine exponent
of $Y$  is $(a',0,\mu_{P_2})$. \hfill\fbox

%\vskip 0.2cm
%\noindent {\bf Proof of Proposition \ref{pro-3.3}.} If $X$ satisfies (H), then $Y$ also satisfies (H) by Lemma  \ref{lem-3.1}. Therefore, the drift coefficient %of $Y$ equals zero by Proposition \ref{pro-3.2} and Lemma \ref{lem-3.4}.\hfill\fbox

\subsection{Converse of Lemma \ref{lem-3.1}}

In this subsection, we consider the converse of Lemma \ref{lem-3.1}. We are particularly interested in the following questions:

\noindent {\bf Question 1.} {\it If for any nonempty proper subspace $A$ of $\mathbf{R}^n$, the projection process $X_A$ of $X$ on $A$  satisfies (H), does  $X$ satisfy (H)?}

\noindent {\bf Question 2.} {\it If  for any one-dimensional subspace $A$ of $\mathbf{R}^n$, the projection process $X_A$ of $X$ on $A$  satisfies (H), does  $X$ satisfy (H)?}

\noindent {\bf Question 3.} {\it If the one-dimensional projection process of $X$ on each coordinate-axis  satisfies (H), does  $X$ satisfy (H)?}

\noindent {\bf Question 4.} {\it Let $S$ be a nonempty proper subspace of $\mathbf{R}^n$. Assume that  the two projection processes $X_S$ and $X_{S^{\perp}}$ of $X$ on $S$ and $S^{\perp}$, respectively, are independent  and satisfy (H). Does  $X$ satisfy (H)?}

\subsubsection{Counterexample for Question 3}

We use a counterexample to show that the answer to Question 3 is
negative.

\begin{exa}\label{exa3.1}
Let $B=(B_t)_{t\geq 0}$ be a standard two-dimensional Brownian
motion. Define
\begin{eqnarray*}
Q=\left(
\begin{array}{cc}
2& 2\\
2& 2
\end{array}\right).
\end{eqnarray*}
Then,
\begin{eqnarray*} \sqrt{Q}=\left(
\begin{array}{cc}
1& 1\\
1& 1
\end{array}\right),
\end{eqnarray*}
and $\sqrt{Q}\mathbf{R}^2=\{(x,y)\in \mathbf{R}^2:x=y\}$.

 Let
$a=(1,-1)^T$ and define $X=(X_t)_{t\geq 0}$ by
$$
X_t=at+\sqrt{Q}B_t,\ \ t\geq 0.
$$
Then, the projection process of $X$ on each coordinate-axis has
non-degenerate Gaussian part and thus satisfies (H) by Theorem
\ref{thm-2.1-5}. On the other hand, the projection process of $X$
on the subspace $\{(x,y)\in \mathbf{R}^2:y=-x\}$ is the uniform
motion, which does not satisfy (H). Therefore, $X$ does not
satisfy (H) by Lemma \ref{lem-3.1}. Note that in this example the
projection process of $X$  on any one-dimensional subspace of
$\mathbf{R}^2$, except for $\{(x,y)\in \mathbf{R}^2:y=-x\}$,
satisfies (H).
\end{exa}

\subsubsection{Partial answer to Question 2}

In this part, we give an affirmative answer to Question 2 under
the assumption that $\mu(\mathbf{R}^n$ $\backslash
\sqrt{Q}\mathbf{R}^n)<\infty$.

\begin{thm}\label{thm-4.7}
 Suppose that $Q$ is degenerate and $\mu(\mathbf{R}^n\backslash \sqrt{Q}\mathbf{R}^n)<\infty$. Then the following three claims are equivalent:

 (i) $X$ satisfies (H);

 (ii) for any one-dimensional subspace $A$ of $\mathbf{R}^n$, the projection process of $X$ on $A$  satisfies (H);

 (iii) the projection process of $X$ on $(\sqrt{Q}\mathbf{R}^n)^{\bot}$ satisfies (H).

\end{thm}
{\bf Proof.} $(i) \Rightarrow (ii)$ and $(iii)$: This is a direct
consequence of  Lemma \ref{lem-3.1}.

$(ii) \Rightarrow (i)$: Let $k$ be the rank of $Q$. Then $0\leq
k<n$.   By the orthogonal transformation of L\'{e}vy processes
(cf. \cite[Section 2.2]{HS12}), we can assume without loss of
generality that $Q={\rm diag}(\lambda_1,\ldots,\lambda_n)$, where
$\lambda_1\geq \lambda_2\geq \cdots\geq
\lambda_k>0,\lambda_{k+1}=\cdots=\lambda_n=0$, and $X$ has the
expression
\begin{eqnarray}\label{thm-4.7-a}
X_t=X_t^{(1)}+X_t^{(2)},\ \ t\geq 0,
\end{eqnarray}
where
\begin{eqnarray*}
&&X^{(1)}_t:=b't+\sqrt{Q}B_t+\int_{\{x\in
\mathbf{R}^k\times\{0\}:|x|\geq 1\}}xN(t,dx)+\int_{\{x\in
\mathbf{R}^k\times\{0\}:|x|< 1\}}x\tilde{N}(t,dx),\\
&&X^{(2)}_t:=\int_{\mathbf{R}^k\times (\mathbf{R}^{n-k}\backslash
\{0\})}xN(t,dx),\nonumber
\end{eqnarray*}
where $b'=-a-\int_{\{|x|<1\}}\mu_1(dx)$, $\mu_1$ is the
restriction of $\mu$ on $\mathbf{R}^n\backslash
\sqrt{Q}\mathbf{R}^n=\mathbf{R}^k\times
(\mathbf{R}^{n-k}\backslash \{0\})$, $N$ is a Poisson random
measure on $\mathbf{R}^+\times (\mathbf{R}^n\backslash\{0\})$
which is independent of the standard Brownian motion
$B=(B_t)_{t\ge 0}$, and $\tilde{N}(t,F)=N(t,F)-t\mu(F)$.

 By the assumption $\mu(\mathbf{R}^n\backslash \sqrt{Q}\mathbf{R}^n)<\infty$, we find that $(X^{(2)}_t)_{t\ge 0}$ is a
 compound Poisson process. Denote $b'=(b'_1,\ldots,b'_k,b'_{k+1},\ldots,b'_n)$ and let ${\cal P}_i$ be the projection of $\mathbf{R}^n$ on the one-dimensional
 subspace $A_i:=\{(x_1,\ldots,x_n)\in \mathbf{R}^n:x_j=0,\ \forall j\neq i\}$ for $i=k+1,\ldots,n$.  Then we
 have that
 $$
 {\cal P}_iX_t=({\cal P}_ib')t+{\cal P}_iX_t^{(2)},
 $$
where $({\cal P}_iX_t^{(2)})_{t\ge 0}$ is a compound Poisson
process and ${\cal P}_ib'\in A_i$.

By (ii), we find that the projection process $({\cal
P}_iX_t)_{t\ge 0}$ satisfies (H) and hence $b'_i=0$ for any
$i=k+1,\dots, n$. Then $b'\in
\mathbf{R}^k\times\{0\}=\sqrt{Q}\mathbf{R}^n$. Therefore, $X$
satisfies (H) by Theorem \ref{thm-2.1-6}.

$(iii)\Rightarrow (i)$: Denote by $P_2$   the projection operator
from $\mathbf{R}^n$ to $(\sqrt{Q}\mathbf{R}^n)^{\perp}$. By
(\ref{thm-4.7-a}), we get
$$
P_2X_t=(P_2b')t+P_2X_t^{(2)},
 $$
where $(P_2X_t^{(2)})_{t\ge 0}$ is a compound Poisson process and
$P_2b'\in (\sqrt{Q}\mathbf{R}^n)^{\bot}$. By (iii), we find that
$(P_2X_t)_{t\ge 0}$ satisfies (H) and hence $P_2b'=0$, which
implies that $b'\in \sqrt{Q}\mathbf{R}^n$. Therefore, $X$
satisfies (H) by Theorem \ref{thm-2.1-6}. \hfill\fbox

\subsubsection{A result on  Question 1}

In this part, we give a partial result on Question 1.

\begin{thm}\label{pro-4.8}
Suppose there exists  a subspace $S$ of $\mathbf{R}^n$ such that $\sqrt{Q}\mathbf{R}^n\subsetneq S \subsetneq \mathbf{R}^n$ and $\mu(\mathbf{R}^n\backslash S)<\infty$.  Then the following three claims are equivalent:

(i) $X$ satisfies (H);

(ii) for any nonempty proper subspace $A$ of $\mathbf{R}^n$, the
projection process of $X$ on $A$  satisfies (H);

(iii) for $A\in \{S,S^\bot\}$,  the projection process of $X$ on  $A$  satisfies (H).
\end{thm}
{\bf Proof.}  By Lemma \ref{lem-3.1}, we get  $(i) \Rightarrow (ii)\Rightarrow (iii)$.

$(iii) \Rightarrow (i)$:  Let $k$ be the dimension of $S$. Then $0<k<n$. By the orthogonal transformation of L\'{e}vy processes, we can assume without loss of
generality that $S=\{(x_1,\dots,x_n)\in \mathbf{R}^n:x_{k+1}=\cdots=x_n=0\}$.  By the assumption that $\sqrt{Q}\mathbf{R}^n\subsetneq S$, we can express $X$ as \begin{eqnarray}\label{thm-4.8-a}
X_t=X_t^{(1)}+X_t^{(2)},\ t\geq 0,
\end{eqnarray}
where
\begin{eqnarray*}
&&X^{(1)}_t:=b't+\sqrt{Q}B_t+\int_{\{x\in
\mathbf{R}^k\times\{0\}:|x|\geq 1\}}xN(t,dx)+\int_{\{x\in
\mathbf{R}^k\times\{0\}:|x|< 1\}}x\tilde{N}(t,dx),\\
&&X^{(2)}_t:=\int_{\mathbf{R}^k\times (\mathbf{R}^{n-k}\backslash
\{0\})}xN(t,dx),\nonumber
\end{eqnarray*}
where $b'=-a-\int_{\{|x|<1\}}\mu_S(dx)$, $\mu_S$ is the
restriction of $\mu$ on $\mathbf{R}^n\backslash
S=\mathbf{R}^k\times
(\mathbf{R}^{n-k}\backslash \{0\})$, $N$ is a Poisson random
measure on $\mathbf{R}^+\times (\mathbf{R}^n\backslash\{0\})$
which is independent of the standard Brownian motion
$B=(B_t)_{t\ge 0}$, and $\tilde{N}(t,F)=N(t,F)-t\mu(F)$.

 By the assumption $\mu(\mathbf{R}^n\backslash S)<\infty$, we find that $(X^{(2)}_t)_{t\ge 0}$ is a
 compound Poisson process. Denote $b'=(b'_1,\ldots,b'_k,b'_{k+1},\ldots,b'_n)$. Let ${\cal P}_S$ and ${\cal P}_{S^\bot}$ be the projections of $\mathbf{R}^n$ on the  subspaces $S$ and $S^\bot$, respectively.  Then we
 have that
 $$
 {\cal P}_{S^\bot}X_t=({\cal P}_{S^\bot}b')t+{\cal P}_{S^\bot}X_t^{(2)},
 $$
where $({\cal P}_{S^\bot}X_t^{(2)})_{t\ge 0}$ is a compound Poisson
process and ${\cal P}_{S^\bot}b'\in {S^\bot}$.

By (iii), we find that the projection process $({\cal
P}_{S^\bot}X_t)_{t\ge 0}$ satisfies (H) and hence ${\cal P}_{S^\bot}b'=0$. Then $b'\in
\mathbf{R}^k\times\{0\}=S$. Thus
 \begin{eqnarray*}\label{thm-4.8-b}
 P_SX^{(1)}_t=X^{(1)}_t,\ \ t\geq 0,
 \end{eqnarray*}
 which implies that
 $$
 P_SX_t=X^{(1)}_t+P_SX^{(2)}_t,\ \ t\geq 0,
  $$
 where $(P_SX^{(2)}_t)_{t\ge 0}$ is a compound Poisson process. By (iii) and Theorem \ref{thm-2.1-16} (iii), we conclude that $(X^{(1)}_t)_{t\ge 0}$ satisfies (H).

 Suppose that $F$ is a semipolar set of $X$. Note that $(X^{(1)}_t)_{t\ge 0}$ satisfies (H),  $(X^{(2)}_t)_{t\ge 0}$ is a  compound Poisson process,  $(X^{(1)}_t)_{t\ge 0}$ and $(X^{(2)}_t)_{t\ge 0}$ are independent. Following the proof of \cite[Theorem 1.2, $(ii) \Rightarrow (i)$]{HS12},  we can show that $F$ is a polar set of $X$ by (\ref{thm-4.8-a}).  Therefore, $X$ satisfies (H). \hfill\fbox

 \begin{rem} In view of Theorem \ref{pro-4.8}, we point out that the additional assumption that the projection process of $X$ on $S^{\bot}$ satisfies (H) should be added to
\cite[Lemma 3.5]{HS18} in order that its conclusion holds. But \cite[Lemma 3.6]{HS18} is still true and thus \cite[Theorem 3.1]{HS18} (i.e., Theorem \ref{thm-2.1-13} of Section 2) holds.
 \end{rem}

\subsubsection{Two propositions on Question 4}

By the orthogonal transformation of L\'{e}vy processes, we find that Question 4 is equivalent to the following question:

\noindent {\bf Question 4'.} {\it Let  $X_1$ and $X_2$ be L\'{e}vy processes on $\mathbf{R}^n$ and $\mathbf{R}^m$, respectively. Suppose that  $X_1$ and $X_2$ are independent and both of them satisfy (H). Does the $\mathbf{R}^{n+m}$-valued  L\'{e}vy process $X=(X_1,X_2)$ satisfy (H)?}

Denote by $\psi_1$ and  $(a_1,Q_1,\mu_1)$ the L\'{e}vy-Khintchine exponent of $X_1$,  by $\psi_2$ and  $(a_2,Q_2,\mu_2)$ the L\'{e}vy-Khintchine exponent of $X_2$, and by $\psi$ and  $(a,Q,\mu)$ the L\'{e}vy-Khintchine exponent of $X$.
Define the $(n+m)\times (n+m)$ matrix
\begin{eqnarray}\label{a}
\bar{Q}:=
\left(
\begin{array}{cc}
Q_1& 0\\
0 & Q_2
\end{array}
\right),
\end{eqnarray}
and two measures $\bar{\mu}_1$ and $\bar{\mu}_2$ on $\mathbf{R}^{n+m}$ by
\begin{eqnarray}\label{b}
&&\bar{\mu}_1(A):=\mu_1(\{x\in \mathbf{R}^n:(x,0_m)\in A\}),\ \ \bar{\mu}_2(A):=\mu_2(\{y\in \mathbf{R}^m:(0_n,y)\in A\}),
\end{eqnarray}
where $0_n$ and $0_m$ are zero elements of $\mathbf{R}^n$ and $\mathbf{R}^m$, respectively, and $A$ is an arbitrary Borel subset of $\mathbf{R}^{n+m}$. By direct calculation, we get
\begin{eqnarray}\label{b-1}
\psi(x,y)=\psi_1(x)+\psi_2(y),\ \ (x,y)\in \mathbf{R}^{n+m},
\end{eqnarray}
and
\begin{eqnarray}\label{c}
a=(a_1,a_2),\ \ Q=\bar{Q},\ \ \mu=\bar{\mu}_1+\bar{\mu_2}.
\end{eqnarray}

\begin{pro} (\cite[Lemma 3.6]{HS18})
If $X_1$ satisfies (H) and $X_2$ is a compound Poisson process, then $X=(X_1,X_2)$ satisfies (H).
\end{pro}

\begin{pro}
If both $X_1$ and $X_2$ satisfy (S), then $X=(X_1,X_2)$ satisfies (H).
\end{pro}
\noindent{\bf Proof.} By the assumption we have that
$$
\mu_1(\mathbf{R}^n\backslash \sqrt{Q_1}\mathbf{R}^n)<\infty,\ \ \mu_2(\mathbf{R}^m\backslash \sqrt{Q_2}\mathbf{R}^m)<\infty,
$$
which together with (\ref{a}), (\ref{b}) and (\ref{c}) implies that
\begin{eqnarray}\label{pro-5-b}
\mu(\mathbf{R}^{n+m}\backslash \sqrt{Q}\mathbf{R}^{n+m})&=&(\bar{\mu}_1+\bar{\mu}_2)(\mathbf{R}^{n+m}\backslash \sqrt{\bar{Q}}\mathbf{R}^{n+m})\nonumber\\
&=&\mu_1(\mathbf{R}^n\backslash \sqrt{Q_1}\mathbf{R}^n)+\mu_2(\mathbf{R}^m\backslash \sqrt{Q_2}\mathbf{R}^m)\nonumber\\
&<&\infty.
\end{eqnarray}
By Theorem \ref{thm-2.1-6}, we find that both $X_1$ and $X_2$ satisfy the Kanda-Forst condition. Then $X$ satisfies the Kanda-Forst condition by (\ref{b-1}). Therefore,  $X=(X_1,X_2)$ satisfies (H) by  Theorem \ref{thm-2.1-6} and (\ref{pro-5-b}).\hfill\fbox

%
%
%
%
%\begin{pro}
%If $X_1$ has bounded resolvent density and $X_2$ satisfies the Kanda-Forst condition, then $X=(X_1,X_2)$ satisfies (H). ???
%\end{pro}

\section{Energy for multidimensional L\'{e}vy processes}\setcounter{equation}{0}

Let $X=(X_t)_{t\geq 0}$ be an $\mathbf{R}^n$-valued L\'{e}vy process on $(\Omega,{\cal F},P)$ with L\'{e}vy-Khintchine exponent $\psi$. For a finite  measure $\nu$ on $\mathbf{R}^n$ and $\lambda>0$, we define its $\lambda$-energy $E_X^{\lambda}(\nu)$ by
$$
E_X^{\lambda}(\nu)=\int_{\mathbf{R}^n} {\rm Re}([\lambda+\psi(z)]^{-1})|\hat\nu(z)|^2dz=\int_{\mathbf{R}^n}\frac{\lambda+{\rm Re}\psi(z)}{|\lambda+\psi(z)|^2}|\hat\nu(z)|^2dz.
$$

Energy plays a fundamental role in the study of Hunt's hypothesis (H). Kanda and Rao gave the following remarkable result on the relation between a measure which does not charge semipolar sets and its energy.

\begin{thm}\label{KR} (Kanda \cite{Knew} and Rao \cite{R88}) Assume that $X$ has resolvent densities with respect to $m_n$. Let $\nu$ be a finite measure which charges no semipolar sets and $E_X^{\lambda}(\nu)<\infty$ for $\lambda>0$.
Then
$$\lim_{\lambda\to\infty}E_X^{\lambda}(\nu)=0.
$$
\end{thm}

Based on Rao \cite{R88}, we proved the following result.

\begin{thm}\label{thm-4.1} (\cite[Theorem 5.1]{HSZ15}) Assume that $X$ has resolvent densities with respect to $m_n$.
Then (H) holds if and only if
$$\lim_{\lambda\to\infty}E_X^{\lambda}(\nu)=0
$$
for any finite measure $\nu$ of finite 1-energy.
\end{thm}

Kanda considered in \cite{Ka91} the space-time process $Y=(Y_t)_{t\geq 0}$ over $X$, which means that $Y$  is a L\'{e}vy process on $\mathbf{R}^1\times \mathbf{R}^n$ defined on the probability space $(\mathbf{R}^1\times \Omega,P^{r,x})$, where $P^{r,x}=\delta_r\otimes P^x$ with $\delta_r$ being the Dirac measure at $r\in \mathbf{R}^1$. The trajectory $Y_t(r,\omega)$ is $(r+t,X_t(\omega))$ and  the L\'{e}vy-Khintchine exponent of $Y$ is $\psi(z)-it$.

\begin{thm} \label{thm-4.3}(\cite[Theorem]{Ka91}) Let $X$ be a L\'{e}vy process on $\mathbf{R}^n$ with transition probability densities and $Y$ be the space-time process over $X$. Let $\nu$ be a finite measure on $\mathbf{R}^1\times \mathbf{R}^n$ of compact support.

(I) Assume that the $\lambda$-energy of $\nu$ for $Y$ is finite. Then,

(i) The $\mathbf{R}^n$-marginal $\nu_2$ of $\nu$ (i.e., $\nu_2(B)=\nu(\mathbf{R}^1\times B)$) has finite $\lambda$-energy for $X$.

(ii) If the $\mathbf{R}^1$-marginal $\nu_1$ of $\nu$ (i.e., $\nu_1(B)=\nu(B\times \mathbf{R}^n)$) is singular to the Lebesgue measure on $\mathbf{R}^1$, then the $\mathbf{R}^n$-marginal $\nu_2$ does not charge any semipolar set.

(II) Consider the case that $\nu$ is of the direct product form $\eta\otimes\xi$.

(i) If $\nu$ has finite $\lambda$-energy for $Y$ and $\xi$ is carried by a semipolar set for $X$, then $\eta$ has a $L^2$-density relative to the Lebesgue measure on $\mathbf{R}^1$.

(ii) If $\xi$ is a finite measure of compact support on $\mathbf{R}^n$ with finite $\lambda$-energy for $X$ and it does not charge any semipolar set for $X$, then we can find a singular measure $\eta$ of compact support so that $\nu=\eta\otimes \xi$ has finite $\lambda$-energy for $Y$.
\end{thm}

Using Theorem \ref{thm-4.3}, Kanda obtained the following characterization of semipolar sets.
\begin{cor} (\cite[Corollary]{Ka91}) Let $X$ be a L\'{e}vy process on $\mathbf{R}^n$ which has transition probability densities. Then a closed set $B$ in $\mathbf{R}^n$ is semipolar if and only if
$$
P^x(X_t\in B\ \mbox{for some}\ t\in A)=0
$$
for every $x\in \mathbf{R}^n$ and every set $A\subset (0,\infty)$ of Lebesgue measure 0.
\end{cor}

In this section, we will follow the idea of Kanda \cite{Ka91} to consider energy for multidimensional L\'{e}vy processes.
From now on till the end of this section, we let $X = (X_t)_{t\geq 0}$ be an $\mathbf{R}^n$-valued L\'{e}vy process and $Y = (Y_t)_{t\geq 0}$ be an $\mathbf{R}^m$-valued L\'{e}vy process. Define $Z_t = (X_t,Y_t)$. Assume that $X$ and $Y$ are independent, $X$, $Y$ and $Z$ have resolvent densities relative to Lebesgue measures on $\mathbf{R}^n,\mathbf{R}^m$ and  $\mathbf{R}^{n+m}$,  respectively.  Denote the L\'{e}vy-Khintchine exponents of $X$ and $Y$ by $\Phi$ and $\Psi$, respectively. Then the exponent of $Z$ is $\Phi(x) + \Psi(y) $ for $x \in \mathbf{R}^n$ and  $ y \in \mathbf{R}^m$.

\subsection{Results}

\begin{pro}\label{pro1}
  Suppose that $\nu$ is a finite measure on  $\mathbf{R}^n\times \mathbf{R}^m$ of compact support and  has finite $\lambda$-energy for $Z$. Then

    (i)   the $\mathbf{R}^n$-marginal $\nu_1$ of $\nu$  has finite $\lambda$-energy for $X$;

    (ii) the $\mathbf{R}^m$-marginal $\nu_2$ of $\nu$  has finite $\lambda$-energy for $Y$.
 \end{pro}

 \begin{pro}\label{pro1-1}
  Suppose that $\nu$ is a finite measure on  $\mathbf{R}^n\times \mathbf{R}^m$ of compact support  with the direct product  form $\eta\otimes \xi$ and has finite $\lambda$-energy for $Z$.

  (i) If $\xi$ is carried by a semipolar set for $Y$, then $\eta$ has a $L^2$-density relative to the Lebesgue measure on $\mathbf{R}^n$;

  (ii) If $\eta$ is carried by a semipolar set for $X$, then $\xi$ has a $L^2$-density relative to the Lebesgue measure on $\mathbf{R}^m$.
 \end{pro}

As a direct consequence of Proposition \ref{pro1-1} and \cite[Corollary of Lemma 2.1]{Ka91}, we obtain the following result.

\begin{cor}\label{cor1}
 Suppose that $\nu$ is a finite measure on  $\mathbf{R}^n\times \mathbf{R}^m$ of compact support  with the direct product form $\eta\otimes \xi$ and has finite $\lambda$-energy for $Z$.

  (i) If  $\eta$ is singular to the Lebesgue measure on $\mathbf{R}^n$, then $\xi$ does not charge any semipolar set;

  (ii) If  $\xi$ is singular to the Lebesgue measure on $\mathbf{R}^m$, then $\eta$ does not charge any semipolar set.
\end{cor}

\begin{pro}\label{pro2}
Assume that $X$ satisfies  (H) and $Y$ satisfies  the Kanda-Forst condition.
   Suppose that $\nu$ is a finite measure on $\mathbf{R}^n\times \mathbf{R}^m$ of compact support with the direct product form $\eta    \otimes\xi$ and has finite $\lambda$-energy for $Z$.    Then,
$$\lim_{\lambda\to\infty}E^{\lambda}_Z(\mu)=0.
$$
\end{pro}

%\begin{rem}
%By \cite[Proposition 5.1]{HSZ15}, we know that if for any finite measure $\mu$ of finite 1-energy relative to $Z$, %$\lim_{\lambda\to\infty}E^{\lambda}_Z(\mu)=0$, then $Z$ satisfies (H).
%\end{rem}

\subsection{Proofs}

\indent {\bf Proof of Proposition \ref{pro1}.}  We only prove (i). The proof of (ii) is similar and we omit it here. We assume without loss of generality that $\nu$ is  a probability measure on $\mathbf{R}^n\times \mathbf{R}^m$.  Then $\nu$ can be disintegrated as
$$
\nu(dx,dy)=\nu_1(dx)\nu_2(x,dy),
$$
 where $\nu_1(dx) = \nu(dx\times\mathbf{R}^m)$ and $\nu_2(x,dy)$ are probability measures on $\mathbf{R}^n$ and $\mathbf{R}^m$, respectively. Set
 $$
 f(x,y)=\widehat{\nu_2(x,\cdot)}(y).
 $$
  Then \begin{equation}\label{mnb}\hat{\nu}(x,y) = \widehat{f(\cdot,y)\nu_1}(x).
  \end{equation}

  By the assumption that the $\lambda$-energy of $\nu$ for $Z$ is finite, we get
  $$
  \int_{\mathbf{R}^m}\left(\int_{\mathbf{R}^n}{\rm Re}\left(\frac{1}{\lambda+\Phi(x)+\Psi(y)}\right)|\hat{\nu}(x,y)|^2dx  \right) dy <\infty.
  $$
  It follows that
  \begin{eqnarray}\label{pro1-a}
  \int_{\mathbf{R}^n}{\rm Re}\left(\frac{1}{\lambda+\Phi(x)+\Psi(y)}\right)|\hat{\nu}(x,y)|^2dx<\infty\ \mbox{for almost all}\ y.
  \end{eqnarray}
 By (\ref{mnb}), we get
\begin{eqnarray}\label{pro1-b}
E^{\lambda}_X(f(\cdot,y)\nu_1) =\int_{\mathbf{R}^n}{\rm Re}\left(\frac{1}{\lambda+\Phi(x)}\right)|\hat{\nu}(x,y)|^2dx.
\end{eqnarray}
For $x\in \mathbf{R}^n$ and $y\in \mathbf{R}^m$, we have
\begin{eqnarray}\label{pro1-c}
&&{\rm Re}\left(\frac{1}{\lambda+\Phi(x)+\Psi(y)}\right)\nonumber\\
&&=\frac{\lambda+{\rm Re}\Phi(x)+{\rm Re}\Psi(y)}{(\lambda+{\rm Re}\Phi(x)+{\rm Re}\Psi(y))^2+({\rm Im}\Phi(x)+{\rm Im}\Psi(y))^2}\nonumber\\
&&\geq \frac{\lambda+{\rm Re}\Phi(x)}{(\lambda+{\rm Re}\Phi(x)+{\rm Re}\Psi(y))^2+({\rm Im}\Phi(x)+{\rm Im}\Psi(y))^2}\nonumber\\
&&\geq \frac{\lambda+{\rm Re}\Phi(x)}{2[(\lambda+{\rm Re}\Phi(x))^2+({\rm Re}\Psi(y))^2+({\rm Im}\Phi(x))^2+({\rm Im}\Psi(y))^2]}\nonumber\\
&&=\frac{\lambda+{\rm Re}\Phi(x)}{2[(\lambda+{\rm Re}\Phi(x))^2+({\rm Im}\Phi(x))^2]\left(1+\frac{({\rm Re}\Psi(y))^2+({\rm Im}\Psi(y))^2}{(\lambda+{\rm Re}\Phi(x))^2+({\rm Im}\Phi(x))^2}\right)}\nonumber\\
&&\geq \frac{1}{2\left(1+\frac{({\rm Re}\Psi(y))^2+({\rm Im}\Psi(y))^2}{\lambda^2}\right)}\cdot \frac{\lambda+{\rm Re}\Phi(x)}{(\lambda+{\rm Re}\Phi(x))^2+({\rm Im}\Phi(x))^2}\nonumber\\
&&=\frac{1}{2\left(1+\frac{({\rm Re}\Psi(y))^2+({\rm Im}\Psi(y))^2}{\lambda^2}\right)}\cdot {\rm Re}\left(\frac{1}{\lambda+\Phi(x)}\right).
\end{eqnarray}
Then, we obtain by (\ref{pro1-a})--(\ref{pro1-c}) that
\begin{eqnarray}\label{d}
E^{\lambda}_X(f(\cdot,y)\nu_1) <\infty\ \mbox{for almost all}\  y.
\end{eqnarray}

We have
\begin{eqnarray}\label{rdr}
|\widehat{f(\cdot,y)\nu_1}(x)|^2 = G_1(x,y) + G_2(x,y),
\end{eqnarray}
 where
\begin{eqnarray*}
 G_1(x,y)&=& |\widehat{{\rm Re}f(\cdot,y)\nu_1}(x)|^2 +|\widehat{{\rm Im}f(\cdot,y)\nu_1}(x)|^2,\\
  G_2(x,y)&=& 2\int_{\mathbf{R}^n} \cos \langle x,z\rangle {\rm Im}f(z,y)\nu_1(dz)\int_{\mathbf{R}^n} \sin \langle x,z\rangle {\rm Re}f(z,y)\nu_1(dz)\\
   &&-2\int_{\mathbf{R}^n} \cos \langle x,z\rangle {\rm Re}f(z,y)\nu_1(dz)\int_{\mathbf{R}^n}  \sin \langle x,z\rangle {\rm Im}f(z,y)\nu_1(dz).
\end{eqnarray*}
It follows that $G_1(x,y) = G_1(-x,y)$, $G_2(x,y) = -G_2(-x,y) $, which  together with ${\rm Re}([\lambda+\Phi(x)]^{-1}) = {\rm Re}([\lambda+ \Phi(-x)]^{-1})$, (\ref{rdr}) and (\ref{d}) implies that for any $\varsigma > 0$,
\begin{eqnarray*}
&&\int_{\{|x|<\varsigma\}}{\rm Re}([\lambda+\Phi(x)]^{-1})G_1(x,y)dx\\
&&=\int_{\{|x|<\varsigma\}}{\rm Re}([\lambda+\Phi(x)]^{-1})(G_1(x,y)+G_2(x,y))dx\\
&&=\int_{\{|x|<\varsigma\}}{\rm Re}([\lambda+\Phi(x)]^{-1})|\widehat{f(\cdot,y)\nu_1}(x)|^2dx\\
&&\leq E^{\lambda}_X(f(\cdot,y)\nu_1)\\
&&<\infty\ \ \mbox{for almost all}\ y.
\end{eqnarray*}
Thus $E_X^{\lambda}({\rm Re}f(\cdot,y)\nu_1) <\infty$  for almost all $y$, where
\begin{equation}\label{asd}
{\rm Re}f(x,y)=\int_{\mathbf{R}^m} \cos\langle y,z\rangle\mu_2(x,dz).
\end{equation}
By (\ref{asd}) and the assumption that $\nu$ has compact support, we find that there exist constants $c>0$ and $\varepsilon>0$
 such that ${\rm Re}f(x,y)>c$  for every $|y|<\varepsilon$ and all $x\in {\mathbf{R}^n}$. Therefore, $E^{\lambda}_X(\nu_1)<\infty$ by \cite[Corollary of Lemma 2.1]{Ka91}.\hfill\fbox

\bigskip
\noindent {\bf Proof of Proposition \ref{pro1-1}.} We only prove (i). By Proposition \ref{pro1}, we get $E_Y^{\lambda}(\xi)<\infty$. If $\xi$ charges a semipolar set, then it charges a compact set $K\subset \mathbf{R}^m$ such that $K\subset \{y\in \mathbf{R}^m: E^y[\exp(-\lambda T_K)]<\delta\}$ for some $\delta\in (0,1)$. Let $\xi_K$ be the restriction of $\xi$ on $K$. By \cite[Corollary of Lemma 2.1]{Ka91}, we have that $E^{\lambda}_Y(\xi_K)\leq E^{\lambda}_Y(\xi)<\infty$. Then $K$ must be non-polar for $Y$ by \cite[Lemma 2.3]{Ka91}. Hence $C^{\lambda}_Y(K)\uparrow C$ as $\lambda\uparrow\infty$ for some positive finite constant $C$ by \cite[Lemma 2.4]{Ka91}, where $C^{\lambda}_Y(K)$ is the {\it $\lambda$-capacity} of $K$ relative to the process $Y$. Thus, we obtain by \cite[Lemma 2.2]{Ka91} that
\begin{eqnarray*}\label{pro1-1-a}
\lim_{\lambda\rightarrow\infty}E_Y^{\lambda}(\xi)\geq \lim_{\lambda\rightarrow\infty}E_Y^{\lambda}(\xi_K)\geq \frac{(2\pi)^m\xi(K)^2}{2C},
\end{eqnarray*}
i.e.,
\begin{eqnarray}\label{pro1-1-a}
\lim_{\lambda\rightarrow\infty}\int_{\mathbf{R}^m}{\rm Re}\left(\frac{1}{\lambda+\Psi(y)}\right)|\hat\xi(y)|^2dy\geq \frac{(2\pi)^m\xi(K)^2}{2C}.
\end{eqnarray}

Similar to (\ref{pro1-c}), we can show that  for $x\in \mathbf{R}^n$ and $y\in \mathbf{R}^m$,
\begin{eqnarray}\label{pro1-1-b}
{\rm Re}\left(\frac{1}{\lambda+\Phi(x)+\Psi(y)}\right)\geq
\frac{1}{2\left(1+\frac{({\rm Re}\Phi(x))^2+({\rm Im}\Phi(x))^2}{\lambda^2}\right)}\cdot {\rm Re}\left(\frac{1}{\lambda+\Psi(y)}\right).
\end{eqnarray}
It follows from (\ref{pro1-1-a}) and (\ref{pro1-1-b}) that
\begin{eqnarray*}
\liminf_{\lambda\rightarrow\infty}\int_{\mathbf{R}^m}{\rm Re}\left(\frac{1}{\lambda+\Phi(x)+\Psi(y)}\right)|\hat\xi(y)|^2dy\geq \frac{(2\pi)^m\xi(K)^2}{4C},\ \ x\in \mathbf{R}^n.
\end{eqnarray*}
By Fatou's lemma, we get
\begin{eqnarray*}
\lim_{\lambda\to\infty}E_Z^{\lambda}(\mu)&=&\lim_{\lambda\to\infty}\int_{\mathbf{R}^n}\left(\int_{\mathbf{R}^m}{\rm Re}\left(\frac{1}{\lambda+\Phi(x)+\Psi(y)}\right)|\hat{\xi}(y)|^2dy\right)|\hat{\eta}(x)|^2dx\\
&\geq& \int_{\mathbf{R}^n}\liminf_{\lambda\rightarrow\infty}\left(\int_{\mathbf{R}^m}{\rm Re}\left(\frac{1}{\lambda+\Phi(x)+\Psi(y)}\right)|\hat{\xi}(y)|^2dy\right)|\hat{\eta}(x)|^2dx\\
&\geq&\frac{(2\pi)^m\xi(K)^2}{4C}\int_{\mathbf{R}^n}|\hat{\eta}(x)|^2dx.
\end{eqnarray*}
Therefore,  $\hat{\eta}\in L^2(\mathbf{R}^n)$, which implies that $\eta$ is absolutely continuous relative to the Lebesgue measure on $\mathbf{R}^n$ and its density belongs to $L^2(\mathbf{R}^n)$.\hfill\fbox

\bigskip

%\begin{pro}\label{pro4}
% Suppose that $\mu$ is a finite measure on $\mathbf{R}^n\times \mathbf{R}^m$ of compact support, has finite energy for $Z$, and has the direct product form $\eta    \otimes\nu$. If  $X$  satisfies (H), then  $\eta$  charges no semipolar sets.
%\end{pro}
%{\bf Proof.} By p. 283 of \cite{BG68}, we know that $\eta$ can be expressed as $\eta=\eta_1+\eta_2+\eta_3$, where $\eta_1$ is carried by a polar Borel set, $\eta_2$ is carried by a semipolar
%Borel set but charges no polar set and $\eta_3$ charges no semipolar set.
%
%By the assumption that $X$ satisfies (H), we know that $\eta_2=0$. By Proposition \ref{pro1}, we know that $\eta$ has finite energy for $X$. It follows from \cite[Corollary of Lemma 2.1]{Ka91} that  $\eta_1$ has finite energy for $X$. Thus the support of $\eta_1$ is non-polar
%and so we have $\eta_1=0$. Hence $\eta=\eta_3$, which charges no semipolar sets. \hfill\fbox
%
%\bigskip

\noindent {\bf Proof of Proposition \ref{pro2}.}  By the assumption that $Y$ satisfies  the Kanda-Forst condition and the proof of \cite[Lemma 4.6]{HS18}, we find that there exist two positive constant $c_1$ and $c_2$  such that for any $y\in \mathbf{R}^m$ and any $\lambda>1$,
\begin{eqnarray}\label{pro2-b}
&&(\lambda+{\rm Re}\Phi(x)+{\rm Re}\Psi(y))^2+({\rm Im}\Phi(x)+{\rm Im}\Psi(y))^2\nonumber\\
&\geq& c_1\left((\lambda+{\rm Re}\Phi(x)+{\rm Re}\Psi(y))^2+({\rm Im}\Phi(x))^2\right),\
\end{eqnarray}
and
\begin{eqnarray}\label{pro2-c}
&&(\lambda+{\rm Re}\Phi(x)+{\rm Re}\Psi(y))^2+({\rm Im}\Phi(x))^2\nonumber\\
&&=(\lambda+{\rm Re}\Phi(x)+{\rm Re}\Psi(y))^2+\left[\left({\rm Im}\Phi(x)+{\rm Im}\Psi(y)\right)-{\rm Im}\Psi(y)\right]^2\nonumber\\
&&\geq c_2(\lambda+{\rm Re}\Phi(x)+{\rm Re}\Psi(y))^2+({\rm Im}\Phi(x)+{\rm Im}\Psi(y))^2.
\end{eqnarray}

By Proposition \ref{pro1}, we know that $\eta$ has finite energy for $X$. Then, by the assumption that $X$ satisfies (H) and Theorem \ref{thm-4.1}, we get
\begin{eqnarray}\label{pro2-d}
\lim_{\lambda\to\infty}E_X^{\lambda}(\eta)=0.
\end{eqnarray}
By (\ref{pro2-b}),   Fatou's lemma (based on (\ref{pro2-c}), the monotonicity of $E^{\lambda}_X(\eta)$ and the assumption that $\nu$ has finite $\lambda$-energy for $Z$),   and (\ref{pro2-d}), we get
\begin{eqnarray*}
&&\lim_{\lambda\to\infty}\int_{\mathbf{R}^{n+m}}{\rm Re}\left(\frac{1}{\lambda+\Phi(x)+\Psi(y)}\right)|\hat\nu(x,y)|^2dxdy\\
&&=\lim_{\lambda\to\infty}\int_{\mathbf{R}^m}\left(\int_{\mathbf{R}^n}\frac{\lambda+{\rm Re}\Phi(x)+{\rm Re}\Psi(y)}{(\lambda+{\rm Re}\Phi(x)+{\rm Re}\Psi(y))^2+({\rm Im}\Phi(x)+{\rm Im}\Psi(y))^2}|\hat\eta(x)|^2dx\right)|\hat\xi(y)|^2dy\\
&&\leq \limsup_{\lambda\to\infty}\frac{1}{c_1}\int_{\mathbf{R}^m}\left(\int_{\mathbf{R}^n}\frac{\lambda+{\rm Re}\Phi(x)+{\rm Re}\Psi(y)}{(\lambda+{\rm Re}\Phi(x)+{\rm Re}\Psi(y))^2+({\rm Im}\Phi(x))^2}|\hat\eta(x))|^2dx\right)|\hat\xi(y)|^2dy\\
&&=\limsup_{\lambda\to\infty}\frac{1}{c_1}\int_{\mathbf{R}^m}E_X^{\lambda+{\rm Re}\Psi(y)}(\eta)|\hat\xi(y)|^2dy\\
&&\le\frac{1}{c_1}\int_{\mathbf{R}^m}\limsup_{\lambda\to\infty}E_X^{\lambda+{\rm Re}\Psi(y)}(\eta)|\hat\xi(y)|^2dy\\
&&=0.
\end{eqnarray*}
The proof is complete.\hfill\fbox

\section{Open questions}

Besides Questions 1, 2, and 4 (equivalently, 4') given in \S 3.2, we would like to present a few more questions on Hunt's hypothesis (H) for further study.

\vskip 0.2cm

\noindent {\bf Question 5.} Does any pure jump subordinator with infinite L\'{e}vy measure satisfy (H)? More generally, does any pure jump $n$-dimensional L\'{e}vy process satisfy (H)? The latter question has been presented in \cite{HSW19}.

\noindent {\bf Question 6.} Dose any one-dimensional L\'{e}vy process with $\int_{\mathbf{R}}(1\wedge |x|)\mu(dx)=\infty$ satisfy (H)?

\noindent {\bf Question 7.} Suppose that $X$ and $Y$ are two independent $n$-dimensional L\'{e}vy processes satisfying (H). Does $X+Y$ satisfy (H)? We have mentioned this question in \cite{HS16}.

\begin{rem}
(i) By Kesten \cite[Theorem 1(f)]{Ke69}, we know that if the answer to Question 7 is yes, then the answer to Question 6 is also yes (see \cite[Section 2.1]{HS18}).

(ii) By Lemma \ref{lem-3.1} and the orthogonal transformation of L\'{e}vy processes,  we know that if the answer to Question 4 is yes, then the answer to Question 7 is also yes.
\end{rem}

\bigskip

{ \noindent {\bf\large Acknowledgments}  This work was supported by National Natural Science Foundation of China (Grant No. 11771309 and No. 11871184), Natural Science and Engineering Research Council of Canada.

\end{document}